**The Babe Ruth Algorithm: a fast, unbiased procedure to randomize presence–absence data matrices with fixed row and column totals**

Giovanni Strona[1], Domenico Nappo[1], Francesco Boccacci[1], Simone Fattorini[2] & Jesus San-Miguel-Ayanz[1]

[1]*Institute for Environment and Sustainability, Joint Research Centre, European Commission, Via E. Fermi 2749, 21027 Ispra (VA), Italy.* [2]*Azorean Biodiversity Group, Departamento de Ciências Agrárias, Universidade dos Açores, Angra do Heroísmo, Terceira, Azores, Portugal*

Correspondence and requests for materials should be addressed to G.S. (email: giovanni.strona@jrc.ec.europa.eu)

A well-known problem in numerical ecology is how to recombine presence-absence matrices without altering row and column totals. A few solutions have been proposed, but all of them present some issues in terms of statistical robustness (i.e. their capability to generate different matrix configurations with the same probability) and their performance (i.e. the computational effort they require to generate a null matrix). Here we introduce the 'Babe Ruth Algorithm', a new procedure that differs from existing ones in that it focuses rather on matrix information content than on matrix structure. We demonstrate that the algorithm can sample uniformly the set of all possible matrix configurations requiring a computational effort orders of magnitude lower than that required by available methods, making it possible to easily randomize matrices larger than $10^8$ cells.

Many ecological patterns are investigated by analyzing presence-absence (1-0) matrices[1]. This holds both for matrices describing distribution of taxonomic units (i.e. species per site matrices) and for matrices describing species interactions (i.e. ecological network matrices)[2]. Species co-occurrence, nestedness, and modularity (i.e. the non-random occurrence of densely connected, non-overlapping species subsets) are typical examples of matrix analyses used in community ecology, biogeography, conservation, host-parasite interaction studies[3]. Although several metrics have been proposed to quantify peculiar matrix aspects, none of them has associated probability levels. Thus, comparing the patterns found in a given matrix with those emerging from randomly generated matrices (null models) is a common procedure to assess the probability that the observed pattern can be explained by chance alone[4]. For this, the matrix under study is randomized several times, and the pattern observed in the original matrix is compared to the patterns observed in the random matrices. In doing this, however, the choice of the 'rules' controlling the randomization process may lead to very different outcomes, thus making the interpretation of results difficult[5].

This happens because null models may vary a lot in their restrictiveness, i.e. in how much a null matrix preserves features of the original matrix[6]. Some basic matrix properties, such as size (i.e. number of matrix cells) and fill (i.e. the ratio between number of occurrences and matrix size) are commonly retained in the generation of null matrices even by the least restrictive null models[5] and several studies highlight the importance of generating null matrices with the same row and column totals as those of the matrix under study to minimize the risk of Type II errors[5,7]. Moreover, preserving row and column totals has deep ecological implications[2]. For example, in a species per site matrix, a row with several presences may suggest that the corresponding species, which occurs in several localities, is an ecological generalist. On the other hand, a column with many presences may indicate that the considered locality, which is inhabited by several different species, has a high resource availability. Thus, preserving row and column totals should be preferable in most situations[2,8,9], and especially in the analysis of species co-occurrence[4].

However, generating random matrices that preserve row and column totals is far from trivial. Several statistical approaches have been proposed to this purpose[10-14].

Most of them use algorithms based on swaps of 'checkerboard units'[13], where a checkerboard unit[15], is a $2 \times 2$ submatrix in one of the two alternative configurations:

$$\begin{bmatrix} 0 & 1 \\ 1 & 0 \end{bmatrix} \text{ and } \begin{bmatrix} 1 & 0 \\ 0 & 1 \end{bmatrix}$$

Swap methods, at each cycle, extract at random two rows and two columns from the matrix. If the $2 \times 2$ submatrix including the cells at the intersection of these rows and columns is a checkerboard unit, it is swapped, i.e. the values of one of its diagonals are replaced with the values of the other one (thus passing from one checkerboard configuration to another one). The most commonly used swap algorithm is the 'sequential swap', which generates a first null matrix by attempting 30000 swaps, and then creates each subsequent null matrix by performing a single swap on the last generated matrix[13]. As a consequence, independently of the number of swaps performed, each matrix in the resulting set differs from the previous one by only 4 matrix elements. Thus, in order to create a set of matrices truly representative of the whole set of possible configurations, a very large number of null matrices must be generated[16,17]. Alternatively, one may use the independent swap algorithm, which creates each null matrix by performing a certain number of swaps (ideally more than 30000) on the original matrix[13,18]. However, swap methods tend to be biased towards the construction of segregated matrices, i.e. matrices with more checkerboard units than expected by chance[14]. Thus none of these approaches (i.e. generating a large number of null matrices, or using the independent swap algorithm) guarantees that the set of null matrices generated constitutes a truly random sample of the universe of possible matrix configurations[13].

A different method, based on random walks on graphs, has been demonstrated to be unbiased[19]. However, this approach, which was primarily conceived as an aid for testing biogeographical hypotheses, has rarely been applied to practical ecological studies, perhaps as a consequence of its complexity. A simpler algorithm is the 'trial-swap'[14], which consists in a subtle modification of the traditional swap algorithms, based on the principle of setting *a priori* the number of total swaps to be attempted on a matrix. The trial-swap algorithm samples different null matrices with the same probability. Thus, the fast algorithms are used to create a first null matrix and subsequent perturbations based on trial swaps are used to generate the uniform distribution of the results[14,20]. Differently, some authors suggested to use the sequential swap and correct it for unequal frequencies[16,19]. More recently, an efficient algorithm for sampling exactly from the uniform distribution over binary or non-negative integer matrices with fixed row and column sums was developed[21]. However, this method can be applied only to small matrices (with a total number of rows and column smaller than 100) or larger matrices with low fill[21]. Thus, despite the higher statistical robustness of some of the above mentioned procedures, the sequential swap still remains the most commonly used algorithm to generate null matrices with fixed row and column totals[17].

Here we show how a childhood pastime may suggest an alternative and effective solution to the problem of presence-absence matrix permutations. Let's imagine that a group of children meets to exchange baseball cards. Let's assume that all cards have the same value, so that a fair trade is one in which one card is given for one card received, and that no one is interested in holding duplicate cards. This scenario may be represented by a matrix with rows corresponding to boys and columns corresponding to cards. Each cell in that matrix will be filled (i.e. equal to 1) if the boy of the respective row owns the card of the respective column, and empty (i.e. equal to 0) if he does not. Now imagine that two children, after comparing their decks of cards, make a trade according to the above rules, i.e. a card is given for one card received, and no trade that leads a boy to own more than one card per type is made. Let's assume that some trades are feasible, i.e. that there are some pairs of children where one of the two owns some cards not owned by the other and vice-versa, and let's have a look at the situation after a few of the possible trades have taken place.

No card has been created or destroyed during the whole trade process. Because of the fair-trade rule, the number of cards owned by each boy has not changed. And because of the one-card-per-type rule, the number of boys owing a particular card has not changed as well. Thus, if we put back boys and cards in rows and column, we will obtain a matrix with a configuration different from the starting one but with the same row and column totals.

On the basis of this idea, we have developed a new algorithm that is computationally not intensive and produces uniformly distributed null matrices with fixed row and column totals. Thanks to its efficiency it can be applied even to very large matrices (>$10^8$ cells) using ordinary personal computers. We named this algorithm after the legendary Babe Ruth, who appears on some of the rarest and most valuable baseball cards of all time.

## Description of the algorithm

Functioning of the Babe Ruth Algorithm is illustrated in Fig. 1. For the sake of clarity, we will refer in the following to a species per site matrix, but the algorithm can be applied to any binary matrix. We call sites *A1*, *A2*, … *An*, and species *Sp1*, *Sp2*, … *Spn*. If the number of species is higher or equal than that of sites, a set of lists of all the species occurring in each site is created (**A**). Then, two lists are extracted at random from this set (**B**). Throughout the text, we will refer to this kind of process using the term 'pair extraction'. Suppose we have extracted sites *A2* and *A3*. The two lists of species (*A2* and *A3*) are compared, in order to identify the set of species present in *A2* but not in *A3* ($A_{2-3}$) and the set of species present in *A3* but not in *A2* ($A_{3-2}$) (**C**). In our example, $A_{2-3}$ has one element, while $A_{3-2}$ has two elements. Thus one species of *A3* extracted at random from the set $A_{3-2}$ is 'traded' with the one species of *A2* belonging to the set *A₂₋₃* (**D**). In general, for any random pair of lists, a certain number of elements exclusive of a list are traded with an equal number of elements exclusive of the other list. The number of trades for each pair of lists will vary randomly from 1 to n, where n is the size of the smaller of the two sets of exclusive elements.

This mechanism makes it possible to randomize the original set of species-site associations without altering the total numbers of species and sites in the newly generated lists. After the steps **B** to **D** have been reiterated for a certain number of times (**E**), the randomized set of species lists is used to recompile the presence-absence matrix (**F**).

If the number of species is lower than that of the sites, the procedure is the same, with the only difference that pair extractions are performed on a set including the lists of all sites where each species occurs.

## Benchmark testing

**Uniform distribution of null matrices.** Ideally the performance of a randomization algorithm should be measured by assessing how extensively it explores the universe of possible matrix configurations that satisfy the requirement of fixed row and column totals. In other words, the algorithm should be neutral, i.e. it should not show preference for particular matrix configurations[13]. To test the robustness of our algorithm towards this issue, we performed two separated tests. First we replicated a simple experiment originally conceived for the same purpose[13,14] which is based on small matrices, and then we performed an experiment with larger matrices.

The first experiment was based on the five matrices shown in Fig. 2. They represent all the possible configurations of a 3 × 3 matrix with row column totals both equal to [1-2-1]. Starting from any of these five matrices, an unbiased algorithm should generate null matrices corresponding to the five different configurations with the same frequency. As already mentioned, the sequential swap algorithm does not satisfy this requirement, and tends to generate more frequently matrices equal to matrix A of Fig.2, due to the fact that the swap rules provide more pathways to configuration A than they do for the other configurations[14].

We used the Babe Ruth Algorithm to randomize matrix A, by reiterating 100 times the procedure of creating 1000 null matrices and recording the frequency of each matrix configuration as more null matrices were progressively generated. Then, we tested if the distribution frequency of the five configurations was significantly different from 1:1:1:1:1. On average, we obtained a distribution frequency of null matrices statistically not different from the desired uniform distribution after the creation of less than 10 null matrices (average p-value of $\chi^2$ tests = 0.62).

The above test should provide reliable information about the robustness of the algorithm towards unequal sampling of matrix configurations. However, since this test is based, in practice, on a single square matrix of very small size, we felt the need to test the robustness of Babe Ruth Algorithm against matrices of different size, shape and fill.

For this, we performed another experiment based on the same principles of the one above, but using a set of 100 larger random matrices, which we created using a procedure aimed at keeping the number of possible alternative matrix configurations low. The procedure works as follows: First, an empty matrix (i.e. a matrix of all zeros) of random size (ranging from $5 \times 5$ to $15 \times 15$) is built, and a random number of checkerboards (ranging from 1 to 5) is selected *a priori*. Then the matrix is filled by randomly extracting one a cell at a time, and assigning it value 1 only if this addition does not make the total number of checkerboards exceed that selected for the matrix. This trial and error process is reiterated until each row and column of the matrix has at least one presence.

For each of these 100 matrices we generated a set of 1000 null matrices using the Babe Ruth Algorithm (with the number of pair extractions conservatively set at 10000). Finally, we used a $\chi^2$ test test to verify if the frequency of the different configurations in each set of null matrices was significantly different from a uniform distribution. The observed frequencies fitted the expected ones for a uniform distribution in all of the 100 sets, with p-values of $\chi^2$ tests ranging from 0.18 to 1 (mean: 0.60, standard error: 0.02).

**Optimal number of pair extractions.** Miklós and Podani[14] performed a simple experiment on two real species per site matrices to assess the minimum number of swaps necessary to sample different matrix configurations with equal probability. The experiment created different sets of null matrices using an increasing number of swaps and computing, for each set, the average of the total number of checkerboards present in each null matrix. The two real matrices were characterized by a very large number of checkerboards. By increasing the number of swaps, the average number of checkerboards in the set of the null matrices decreased progressively until becoming stable around a certain value. The authors identified the minimum number of swaps needed to reach the uniform distribution in null matrices as the number of swaps necessary to reach a stable value of average checkerboards.

We used the same approach to identify the optimal number of pair extractions to be used with the Babe Ruth Algorithm. First, we applied the Babe Ruth Algorithm to the same two binary matrices used in the original experiment[14]. The first matrix included data for the avifauna of the Vanuatu Archipelago (56 species on 28 islands) [13], while the second included presence-absence of 118 plant species in 80 quadrats of 3 × 3 m located in the Sashegy Nature Reserve, Budapest[22,23]. We computed the total number of checkerboards of any species-site matrix by summing up the number of checkerboard units (*CU*) computed for each possible pair of rows. For each pair of species (i.e. rows), a *CU* value was computed as $(R_i-S) \times (R_j-S)$, where $R_i$ is the total number of occurrences of the *i*-th species, $R_j$ is the total number of occurrences of the *j*-th species, and $S$ is the number of shared sites, i.e. the number of sites where both species occur[15].

For both Vanuatu and Sashegy matrix we generated 100 sets of 10000 null matrices using an increasing number of pair extractions (using an arithmetic progression from 0 to 10000, with common difference of 100). In both experiments the average numbers of checkerboards converged to those reported in the original experiment[14], namely 14060 for the Vanuatu matrix and 549626 for the Sashegy matrix. Results for the two matrices are reported, respectively, in Figure 3A and 3B. Very few pair extractions (less than 1000 for both matrices) were enough to reach the stable value of average checkerboards, corresponding to the uniform distribution of null matrix configurations. For any set of null matrices generated using more than this number of pair extractions, the abundance of checkerboards resulted significantly higher than chance in both matrices, with p<0.0001. The Vanuatu matrix has been already investigated in several papers, with different outcomes[10,12,13]. Our results are consistent with those obtained using methods proven to provide unbiased p-values[14,16,19].

However, the optimal number of pair extractions clearly depends on the size of a matrix, or, better, on its smallest dimension (i.e. the minimum between the number of rows and the number of columns). Thus, using the same approach as above, we estimated the number of pair extractions necessary to reach the stable value of average checkerboards in a large set of real matrices of various sizes. For this we used all the 295 matrices provided together with the Nestedness Temperature Calculator software[24]. For each of these matrices, we generated different sets of 1000 null matrices by using the Babe Ruth Algorithm with an increasing number of pair extractions (using an arithmetic progression starting from a value equal to the smallest dimension of the matrix, with a common difference of 1), until the number of expected checkerboards stabilized, i.e. it did not change by more than 1% in 100 subsequent sets of null matrices. Finally, we compared the size, the minor and major dimensions and the fill of each matrix with the minimum number of pair extractions necessary to reach the stability of average checkerboards. Among all the investigated matrices (n=295), this number was very small, ranging from 3 to 366 (mean: 24.34, standard error: 2.22). Moreover, it was in most cases smaller than or equal to the largest matrix dimension (with an average of 1.2 times the largest matrix dimension). Thus, the value of

1000 pair extractions suggested by the experiment performed on the Vanuatu and the Sashegi matrices is likely to be highly conservative.

**Computational demand of the Babe Ruth Algorithm.** We performed a simple test to compare the difference in terms of computational effort between the Babe Ruth Algorithm and the swap methods. For this we compared the number of pair extractions and the number of swaps necessary to maximally perturb a matrix, i.e. to reach a situation where the percentage of matrix cells in a different position with respect to the starting configuration remains stable when new swaps/pair extraction are performed. To this purpose we used a moderately large random matrix of 100 rows × 100 columns that we created by filling cells with 0 and 1 with the same probability. Then we conducted two separate analyses by performing, respectively, one million pair extractions and one million (dual) swaps on that matrix. In both analyses, we measured at each step the percentage of cells differing from the corresponding ones in the original matrix. After a certain number of pair extractions and swaps, the average percentage of cells with different values from the original configuration became stable (with a value around 50%, Fig. 4A). To reach this plateau, about 50000 swaps were necessary, which is consistent with the number recommended by recent work[17,25]. Crucially, the same value was reached by the Babe Ruth Algorithm after less than 200 iterations (Fig. 4B).

This does not necessary imply that the difference between the two algorithms, in terms of computational effort, scales accordingly, as comparing and recombining the species composition of two areas is more complex than evaluating and swapping a 2 × 2 matrix. Moreover, the amount of data processed by a swap is constant, while that processed by a pair extraction is not. To investigate this aspect, we replicated the above test by recording the time necessary to reach the maximum matrix perturbation using, alternatively, swaps and pair extractions. We recorded only the time necessary to perform single swaps (i.e. the random selection and, if possible, the swap of a 2 × 2 submatrix), and single pair extractions (i.e. the random extraction of two species-area lists, the comparison of the two lists and, if possible, the random exchange of some not-shared elements). We performed this test on a quad-core (Intel Xeon E5-2630 @ 2.30GHz) workstation hosting a 64 bit Linux OS. The Python code necessary to replicate the experiment is reported in the Appendix. Using swaps, the time required to reach the maximum perturbation is around 1 second (Fig. 4C). Using pair extractions, the maximum perturbation is reached about twenty times faster (~0.05 seconds, see Fig. 4D).

Another important aspect to be considered is matrix fill. Most real ecological matrices are characterized by relatively high abundance of zeros. For example, the average fill (± SE) in Atmar and Patterson's set is 0.39±0.01. This means that, on average, more than 60% of the cells in a matrix of this set are empty. Swap algorithms do not take this into account, resulting in a large number of unsuccessful swaps tried on empty matrices. By contrast, the Babe Ruth Algorithm is much more efficient in this sense, as it focuses only on matrix presences, being therefore less affected by matrix size than swap methods. In other words, the differences in terms of computation efficiency between the Babe Ruth Algorithm and the swap methods will be more pronounced as matrix fill decreases. We demonstrated this empirically, by replicating the above experiment on 6 different random matrices of 100 rows × 100 columns with different fill values (respectively 10%, 30%, 50%, 70%, 90%). Matrix fill significantly affected the number of swaps necessary to maximally perturb the matrix. In particular, a large number of swaps (~100000) was necessary for the least and the most filled matrices (Fig. 5A). By contrast, as expected, the differences in matrix fill did not affect the performance of Babe Ruth Algorithm (Fig. 5B). Analogous patterns were observable when computational time (measured as described above) was compared to matrix perturbation (Fig. 5C-D). Thus, pair extractions made it possible to reach the maximum perturbation degree of both the most (90%) and the least (10%) filled matrix about 70 times faster than swaps did. Noteworthy, in all cases, the maximum percentage of matrix perturbation reached using pair extractions was always higher than that reached using swaps. This provides indirect evidence that swap methods do not explore uniformly the universe of possible matrix configurations. However, the difference in performance between the Babe Ruth Algorithm and swap methods is not surprising. Let us imagine a matrix of *R* rows and *C* columns having only one checkerboard. A matrix of this kind admits only one alternative configuration. The probability to find this only configuration by extracting 2 × 2 submatrices is equal to $1/(R^2 \times C^2)$. In contrast, the probability to find it by using pair extractions is $1/R^2$ (if *R*<*C*) or $1/C^2$ (if *R*>*C*).

In general the Babe Ruth Algorithm perturbs the matrix in most of the iterations attempted, while the swap methods have to check several submatrices in order to obtain the same effect. Fig. 6 shows the relationship between the number of attempted swaps/pair extractions and the corresponding number of successful ones (i.e. those producing some change on the matrix) for matrices of various fill. Independently of matrix fill, each pair extraction produces some modification to the matrix. By contrast, for a 50% filled matrix, ~ 8 random attempts are necessary to successfully perform a swap. This number increases symmetrically for less and more filled matrices (~10 for a 30% and a 70% filled matrix, and ~60 for a 10% and a 90% filled matrix).

## Discussion

The statistical challenge that originated this study has interested ecologists for more than 30 years[4]. Exploring the universe of different configurations of matrices with fixed row and column sums has several applications in other scientific fields too, including neurophysiology (in the study of multivariate binary time series), sociology (in the analysis of affiliation matrices), and psychometrics (in item response theory) [21].

Yet, most of the solutions that have been proposed so far present major drawbacks. Unbiased algorithms are computationally highly intensive and therefore slow, while faster algorithms do not ensure that different null matrices are generated with equal frequencies, which may affect result robustness[14].

Another limitation is that most available methods cannot be easily used for very large matrices. What makes these algorithms slow, is that they require a huge number of iterations to properly randomize the original matrix. A recent study demonstrated that more than 50000 swaps are necessary to prevent the occurrence of type I errors in a matrix of average size, and that, for large matrices, many more swaps are needed[17]. The same authors also suggest that, if the dataset under examination is particularly large, the analysis should be repeated by progressively increasing the number of swaps until the p-value stabilizes. In common ecological analyses, a standard number of null matrices to assess the significance of a given pattern is 1000, and in macroecological studies, it is quite common to investigate patterns in tens (sometimes hundreds) of large matrices[26]. Thus, in order to evaluate the significance of ecological patterns in a set of 100 matrices, a total of $5 \times 10^9$ swaps is needed. This is quite a large number, even for modern calculators, so that the task would be inevitably time consuming. However, the main issue connected to these algorithms is that they do not ensure robust results for large matrices unless the number of swaps is raised to a number that makes computation itself rather impractical. To provide support to their observations, Fayle and Manica[17] performed some tests on artificial matrices of 900 rows × 400 columns in order to simulate the maximum size of matrices ever used in a null model analysis. A matrix of this size is arguably large from an ecological perspective. However, the increasing data availability for both geographical distribution and ecological networks offers new possibilities to investigate ecological patterns at very large scale and/or at a very high resolution[27]. This makes it likely that the use of matrices even larger than the above mentioned ones will become quite common in the near future. For this, better performing algorithms are needed.

In this paper we propose a solution to a long debated methodological puzzle. All the tests we performed suggest that the Babe Ruth Algorithm is more efficient than most available tools. Using a simple R function (which is provided as Supplementary Information), the Babe Ruth Algorithm was able to randomize a $10^4 \times 10^4$ matrix in less than 10 seconds on a budget notebook equipped with a dual core processor (Intel Core i3-370M @ 2.40 GHz) running a 64 bit Linux OS. The same task required almost 30 minutes when using the independent swap method (with the optimized function *randomizeMatrix* of the package 'picante'[28]) and setting the number of swaps to the recommended value of twice the number of presences in the matrix multiplied by the number of attempts necessary to successfully perform a swap[14] (since the matrix was approximately 50% filled, we used our estimate of 8 attempts per successful swap).

Despite the fact that the logic behind the Babe Ruth Algorithm and its implementation are quite far from swap methods, the two approaches are in a certain way closely related. In practice, what the Babe Ruth Algorithm does, is (1) to identify in a single step all the possible swaps between a pair of matrix rows (or columns) and (2) to perform any of them with a random probability. The first aspect makes the algorithm efficient, while the second makes it robust towards unequal sampling of random matrices. Although we are convinced ours will not be the last word on the subject, we hope that the principles of our approach will be useful for future research, orienting ecologists' efforts more on the information the data convey, rather than on the way they are distributed on a matrix.

## Acknowledgements

The views expressed are purely those of the writers and may not in any circumstances be regarded as stating an official position of the European Commission. We thank Pieter S.A. Beck for his useful comments on the manuscript.

## Author contributions

G.S. conceived the idea, performed the analyses, developed the code and wrote the manuscript, D.N. and F.B. optimized the code, J.S.M.A. and S.F. reviewed the manuscript, J.S.M.A. coordinated the research team. All the authors participated in discussions of the research.

## Figure legends

**Figure 1** Functioning of the Babe Ruth Algorithm. Given a species per site binary matrix, a set including the lists of all species (columns) occurring in each site (rows) is created (**A**). Two lists are extracted at random from this set (**B**). The two lists are compared, in order to identify the set of species occurring in one list but not in the other and vice versa ($A_{2-3}$ and $A_{3-2}$) (**C**). $A_{2-3}$ has one element, while $A_{3-2}$ has two elements. Thus one species of $A3$ extracted at random from the set $A_{3-2}$ is 'traded' with the one species of $A2$ belonging to the set $A_{2-3}$ (**D**). After the steps **B** to **D** have been reiterated for a certain number of times (**E**), the randomized set of species per site lists is used to recompile the presence-absence matrix (**F**).

**Figure 2** The five possible configurations of a presence-absence matrix with both row and column totals equal to [1, 2, 1].

**Figure 3** Relationship between average number of checkerboard units in sets of 1000 randomizations

of the Vanuatu matrix (**A**) and the Sashegi matrix (**B**), and the corresponding number of pair extractions applied to the Babe Ruth Algorithm to generate each set of random matrices. Grey horizontal lines indicate the expected average numbers of checkerboards as reported by Miklós & Podani (2004), namely 14060 for the Vanuatu matrix and 549626 for the Sashegy matrix.

**Figure 4** Comparisons between computational demand of swaps and pair extractions necessary to perturb a 100 × 100 random matrix. Perturbation is measured as the percentage of cells differing from the corresponding ones of the original matrix. **A**: Relationship between number of performed pair extractions (continuous line) and swaps (dotted line), and the corresponding matrix perturbation degree; **B**: relationship between number of performed pair extractions and the corresponding matrix perturbation degree; **C**: relationships between the computational time (in seconds) of pair extractions (continuous line) and swaps (dotted line), and the corresponding matrix perturbation degree; **D**: relationships between the computational time (in seconds) of pair extractions and the corresponding matrix perturbation degree.

**Figure 5** Comparisons between computational demand of pair extractions and swaps required to perturb five 100 × 100 random matrices with different fill values (i.e. percentage of occupied cells), and the corresponding matrix perturbation degree. Perturbation is measured as the percentage of cells differing from the corresponding ones of the original matrix. **A**: Relationship between number of performed pair extractions (continuous line) and swaps (dotted line), and the corresponding matrix perturbation degree; **B**: relationship between number of performed pair extractions and the corresponding matrix perturbation degree; **C**: relationships between the computational time (in seconds) of pair extractions (continuous lines) and swaps (dotted lines), and the corresponding matrix perturbation degree; **D**: relationships between the computational time (in seconds) of pair extractions

and the corresponding matrix perturbation degree. Colors correspond to different fill values (magenta = 10%; blue = 30%; green = 50%; orange = 70%; red = 90%).

**Figure 6** Relationship between the number of attempted pair extractions (continuous lines) and swaps (dotted lines) and the corresponding number of successful ones (i.e. those producing some change on the matrix) for five 100 × 100 random matrices with different fill values (i.e. percentage of occupied cells). All the 5 lines reporting the relationships for pair extractions are almost perfectly overlapping. A nearly perfect overlap is also observed in the lines reporting the relationships for swaps for the two pairs of matrices with respective fill equal to 10% and 70%, and 30% and 50%.

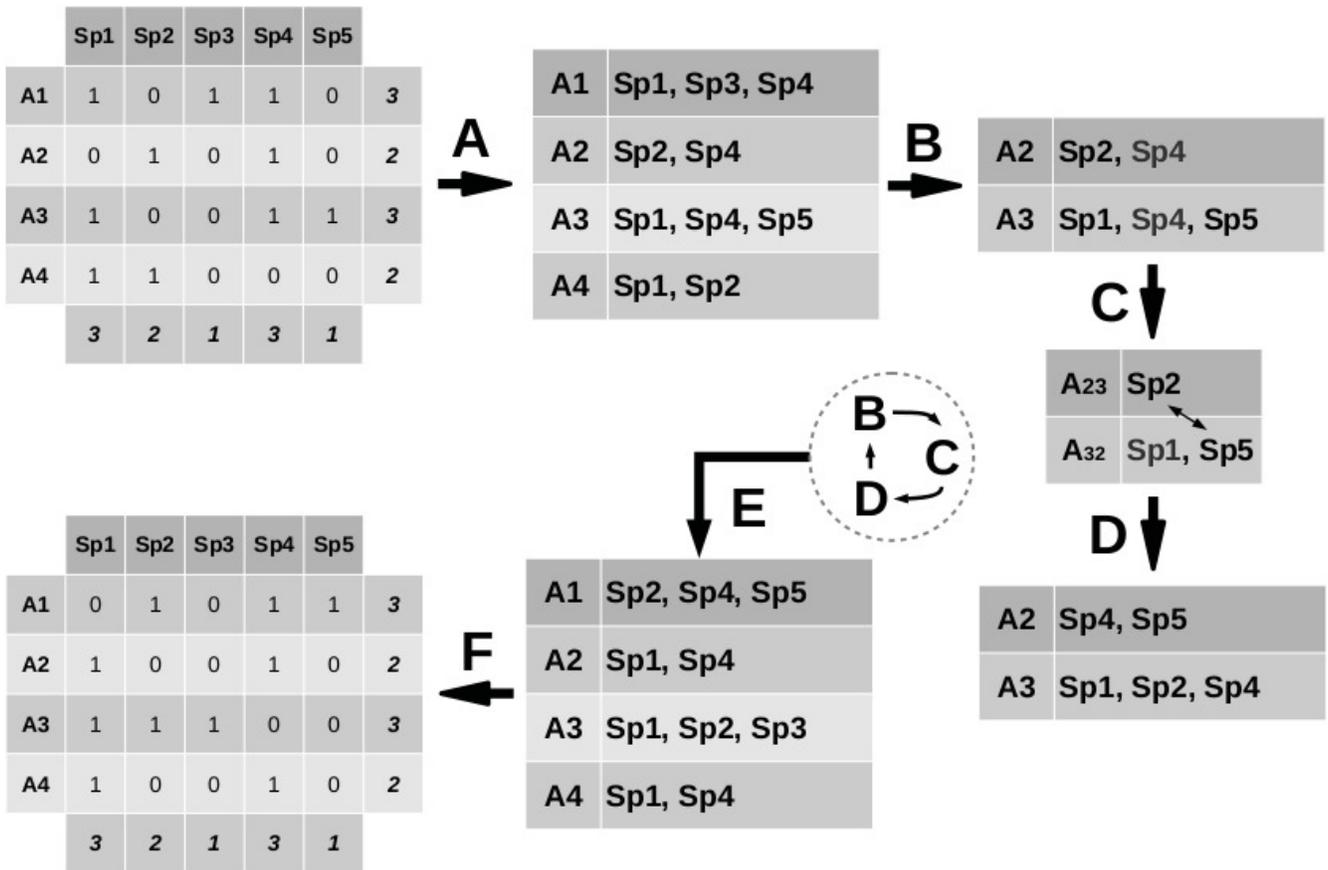

Figure 1

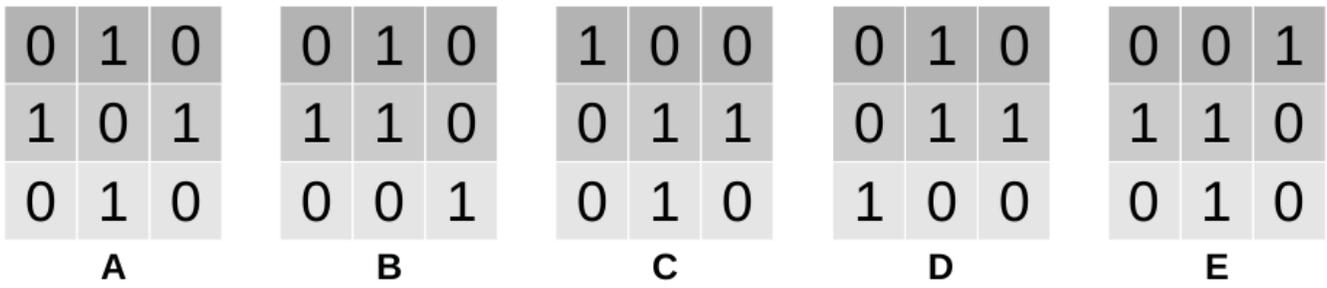

**Figure 2**

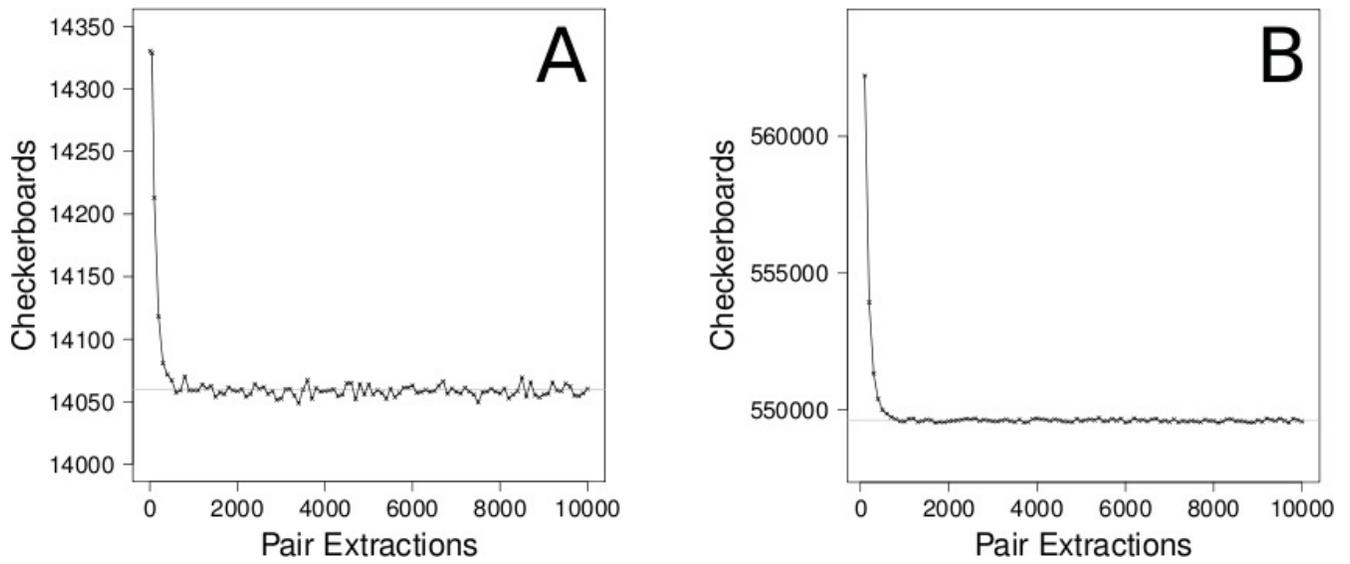

**Figure 3**

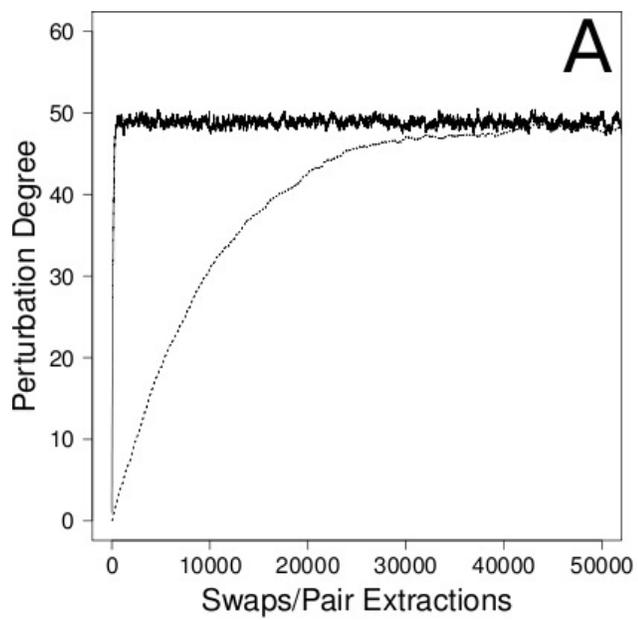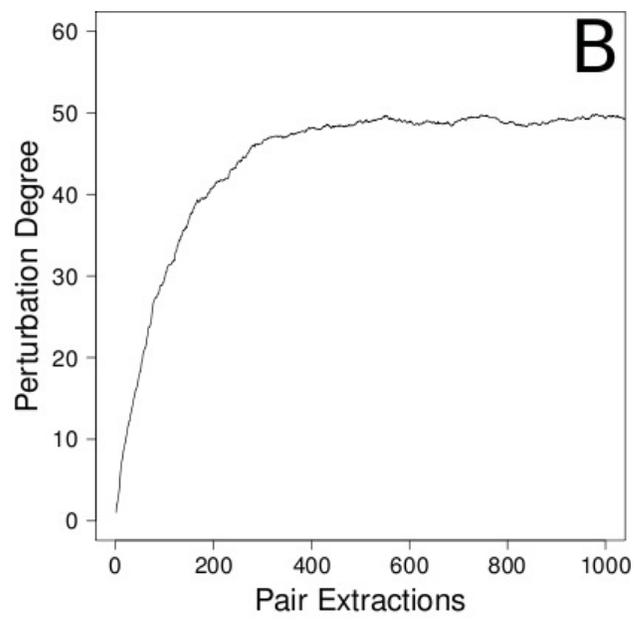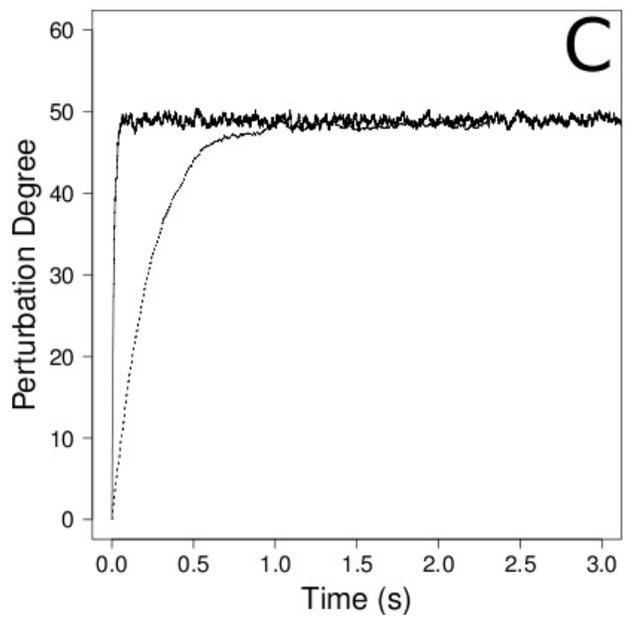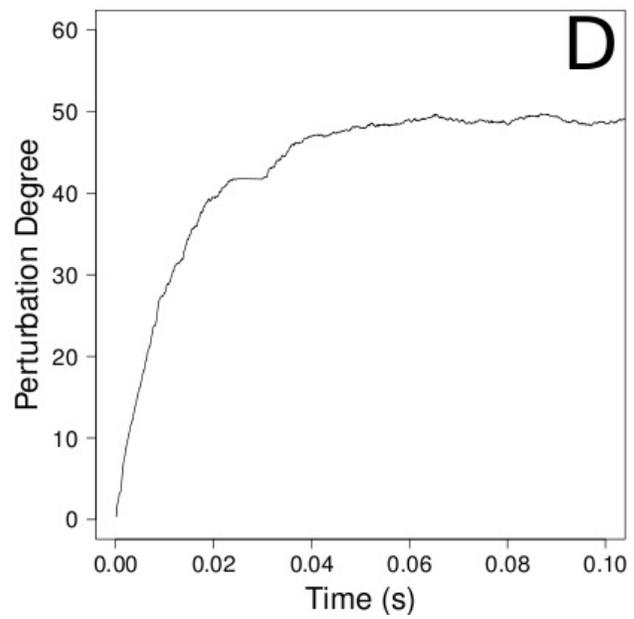

**Figure 4**

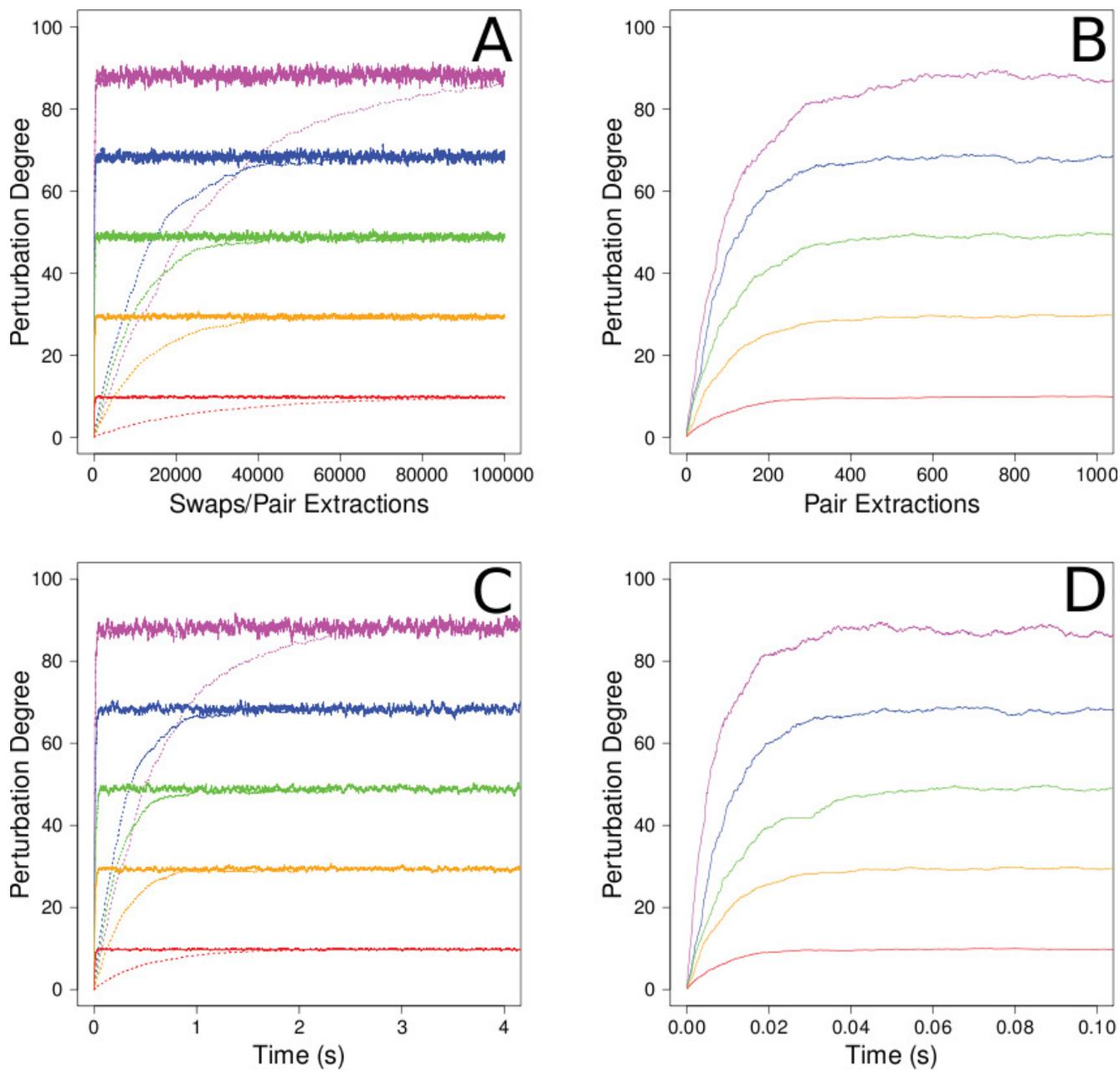

**Figure 5**

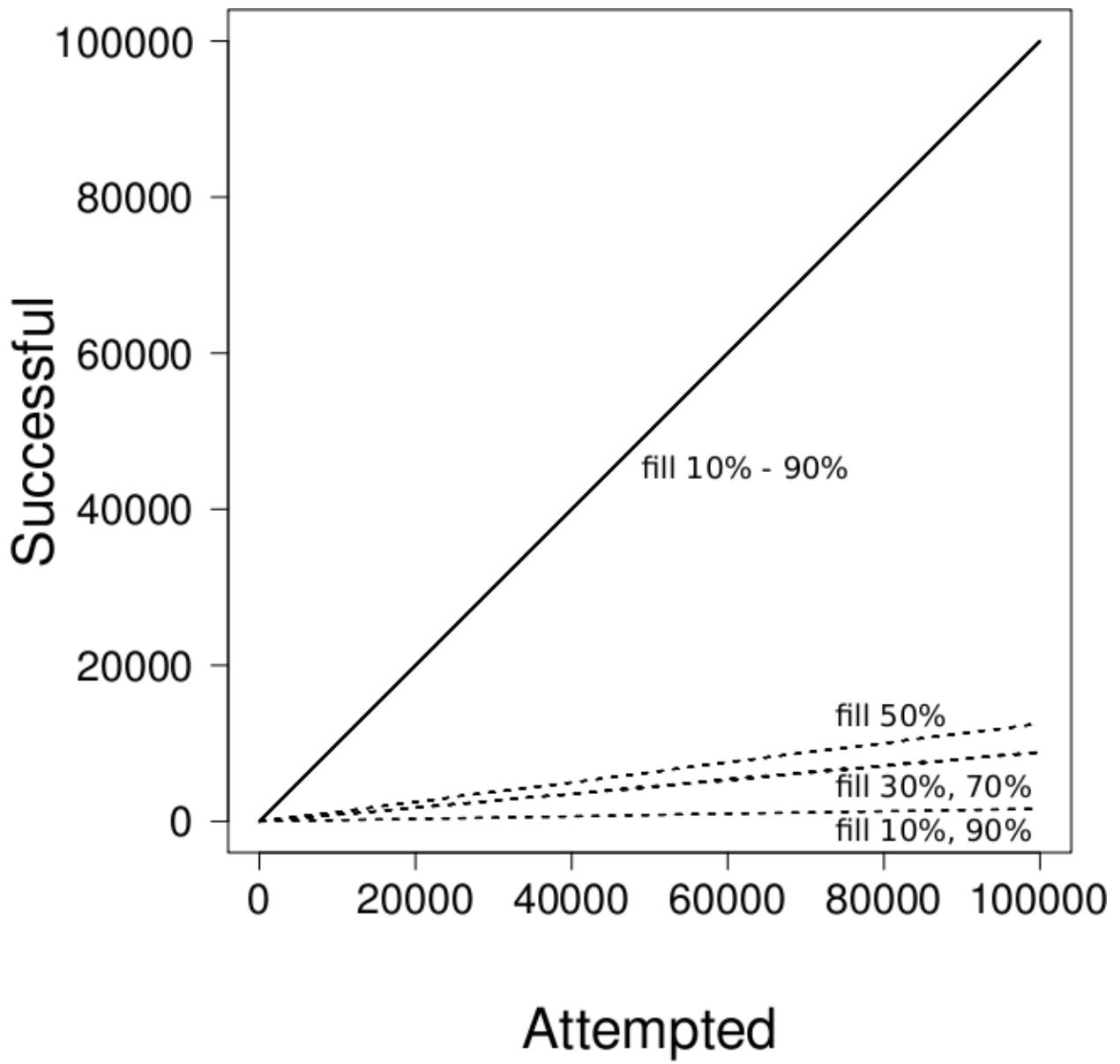

Figure 6